\documentstyle[12pt]{article}
\newtheorem{thm}{Theorem}
\newtheorem{prop}{Proposition}
\newtheorem{cor}{Corollary}
\newtheorem{defn}{Definition}

\newcommand{\T}{Teichm\"{u}ller }
\newcommand{\M}{M\"{o}bius }

\newcommand{\C}{$\hat {\Bbb C}$}
\newcommand{\MM} {PSL($2,\Bbb C$)}

\newcommand{\Bbb}{\bf} 
\newcommand{\sn}{\smallskip\noindent}

\newenvironment{pf}{\noindent {\it Proof.~}}{\hfill$\Box$}
\batchmode
\begin{document}

\title{Coordinates for Teichm\"{u}ller spaces of b-groups with torsion.}
\author{Pablo Ar\'{e}s Gastesi\thanks{AMS Subj. Class. 1991 Primary
32G15, 30F40}}
\date{\today}
\maketitle

\begin{abstract}
In this paper we will use b-groups to construct coordinates for the \T
spaces of $2$-orbifolds. The main technical
tool is the parametrization of triangle groups, which allows us to
compute explicitly formul\ae\ for generators of b-groups uniformizing
orbifolds. In this way, we obtain a technique to pass from
the abstract objects of deformation spaces to concrete
calculations with M\"obius transformations. We explore this computational
character of our work by finding the expressions of certain classical
isomorphisms between \T spaces.
\end{abstract}

One of the most active lines of research in geometric function theory
nowadays deals with the problem of finding parametrizations of \T spaces that
are useful for computations (\cite{ruben:coord}, \cite{kra:horoc},
\cite{mar:coor}, \cite{mas:param1}, \cite{mas:param2}, \cite{min:teic},
\cite{sep:param}). This paper makes a two fold contribution to this
aspect of Complex Analysis: we will give coordinates for the \T spaces of
b-groups with torsion (or equivalently, for the \T spaces of
$2$-orbifolds) and we will use our coordinates to compute explicit
formul\ae\ of well known isomorphisms between deformation spaces. Our
coordinates are good from a computational point of view because,
given a point $\alpha$ in the \T space of a b-group,
$T(p,n;\nu_1,\ldots,\nu_n)$, we have
a technique to compute a set of \M transformations that generate a
Kleinian group $\Gamma$, whose coordinates in the space
$T(p,n;\nu_1,\ldots,\nu_n)$ are $\alpha$.

Given an orbifold $S$, defined over a surface of genus
$p$ with $n$ special points, and a maximal partition $\cal C$ on $S$
(that is, a way of decomposing $S$ into \lq pairs of pants\rq), we can find a
Kleinian group $\Gamma$, acting discontinuously on a simple connected
open set $\Delta$ of $\hat{\Bbb C}$, such that $\Delta/\Gamma\cong
S$. Besides $S$, the group $\Gamma$ uniformizes a finite
number of rigid orbifolds. Therefore, the \T space of $\Gamma$
$T(\Gamma)$ is a model for
$T(S)$, the \T space of the orbifold $S$. This last set is important
in the study of Riemann surfaces because it is the universal covering
space of the Riemann space $R(S)$, which parametrizes the
biholomorphic classes of complex structures on $S$. The advantage of
using the deformation space of the group over the deformation space
of the orbifold lies in the fact that one can do explicit
computations with \M tranformations, obtaining in this way some
properties of $T(S)$.

Using the partition $\cal C$, we can decompose the group $\Gamma$
into a set of subgroups, $\Gamma_1,\ldots,\Gamma_{3p-3+n}$,
with the property that $T(\Gamma_j)$ has dimension $1$.
By a theorem of B. Maskit, the restriction mapping
$T(\Gamma)\rightarrow\prod_{j=1}^{3p-3+n} T(\Gamma_j)$ is one-to-one
and holomorphic. Therefore, to give coordinates on $T(\Gamma)$ it
suffices to study the one-dimensional cases. This is done in detail
in \S $3$ of this paper. Putting together our computations with
the above embedding, we get the following general result, which
generalizes the torsion-free case studied by I. Kra.

\setcounter{thm}{9}
\begin{thm}[\cite{kra:horoc} and \S $3.7$]
Let $\cal S$ be an orbifold with hyperbolic
signature $\sigma=(p,n;\nu_1,\ldots,\nu_n)$,
and let $\cal C$ be a maximal partition on $\cal S$.
Then there exists
a set of global coordinates for the deformation space
$T(S)=T(p,n;\nu_1,\ldots,\nu_n)$, say $(\alpha_1,\ldots,\alpha_d)$,
where $d=3p-3+n$, and a set of non-negative numbers,
$y_1^1,\ldots,y_1^d,y_2^1,\ldots,y_2^d$,
which depends only
on the signature $\sigma$ and the partition $\cal C$, such that
$$\{(\alpha_1,\ldots,\alpha_d)\in {\bf C}^d;~
{\rm Im}(\alpha_i)> y_1^i\; ,\forall \;1\leq i\leq d\}\subset
T(p,n;\nu_1,\ldots,\nu_n)$$
and
$$ T(p,n;\nu_1,\ldots,\nu_n)\subset
\{(\alpha_1,\ldots,\alpha_d)\in {\bf C}^d;~
{\rm Im}(\alpha_i)> y_2^i\; ,\forall \;1\leq i\leq d\}.$$
Given a point $\alpha=(\alpha_1,\ldots,\alpha_d)$ in
$T(p,n;\nu_1,\ldots,\nu_n)$, it is possible to find explictly
a set of $2p+n$ \M transformations that generate a group $\Gamma$,
whose coordinates in that \T space are $\alpha$.\end{thm}

The Kleinian groups that we will use are known as terminal regular
b-groups (see \S $1$ for the definition). The are built from
triangle groups by
a finite number of applications of the Klein-Maskit Combination
Theorems. Therefore, to prove the above theorem, we first need a way of
computing generators for triangle groups. This is the
content of the following result, which is the main technical tool
of our work.\vspace{5mm}

\noindent {\bf Theorem (\S $2.2$ and \S $2.6$)}
{\it Given three distinct points $(a,b,c)$ in the Riemann sphere,
and a signature $\sigma=(0,3;\nu_1,\nu_2,\nu_3)$, which is either hyperbolic
or $(0,3;\infty,2,2)$, there exists a pair of M\"{o}bius
transformation $A$ and $B$, uniquely determined by the parameters
$(a,b,c)$, such that the group generated by them is a triangle group
of signature $\sigma$.
The transformations $A$ and $B$ are explicitly
computable from $(a,b,c)$ and $\sigma$.}\vspace{.5cm}

The explicit character of our coordinates allows us to compute some
classical isomorphisms between different \T spaces as follows.
Let $S$ be a surface of genus
$2$, unformized by a b-group $\Gamma$. Since all surfaces of genus $2$
are hyperelliptic, we have a conformal involution $j$ on $S$ with $6$
fixed points. The quotient orbifold $S'=S/\!<j>$ has signature
$(0,6;2,\ldots,2)$. It is a classical result (\cite{pat:dist})
that the spaces $T(S)$ and
$T(S')$ are biholomorphically equivalent. To find a mapping between them,
we first find the lifting of $j$, say $A_2^{1/2}$,
to the covering determined by $\Gamma$. We have that the group generated by
$\Gamma$ and $A_2^{1/2}$ is a b-group uniformizing $S'$.
Using explicit calculations on $T(\Gamma)$ and $T(\Gamma')$,
we get the following result.\vspace{5mm}

\begin{thm}[\S $4.1$]
The mapping
$$(\tau_1,\tau_2,\tau_3)\mapsto (\frac{\tau_1}{2},1+\tau_2,
1+\frac{\tau_3}{2})$$
gives an isomorphism between $T(2,0)$ and $T(0,6;2,2,2,2,2,2)$.\end{thm}

This paper is organized as follows. \S $1$ contains non-standard
background on Kleinian groups and \T spaces. In \S $2$ we compute
generators for triangle groups.  These computations are used in
\S $3$ to develop coordinates
for \T spaces of terminal regular b-groups; in particular we prove
theorem $10$. In \S $4$ we prove theorem $11$, and indicate how our
methods can be used to prove similar results.

The content of this paper is part of the author's Ph.D. thesis.
I would like to thank my advisor, Irwin Kra, for all his help during my
years as a graduate student; and thank Chaohui Zhang
and Suresh Ggovindarajan for many useful conversations and comments
on my work. I also want to thank the referee for making many
useful comments which have improved a first, and very hard to
understand, version of this paper.

\section{Background}
\noindent {\bf 1.1.} Throughout this paper, we will identify the group of \M
transformations with the projective special linear group,
\MM. The mapping
$z\mapsto{\load{\large}{\displaystyle}\frac{az+b}{cz+d}}$
will be identified with $\left[\begin{array}{cc}a&b\\c&d\end{array}\right]$.
The square brackets denote an element of PSL($2,\Bbb
C$); that is, a class of matrices of SL($2,\Bbb C$). If we take a
particular lifting of an element of PSL($2,\Bbb C$) to SL($2,\Bbb
C$), we will use parenthesis.

An elliptic transformation of finite order $q$ is called {\bf geometric}
if it can be
conjugated in PSL($2,\Bbb C$) to $z\mapsto exp(\pm 2\pi i/q)z$.
Observe that only
\lq minimal rotations\rq\ are geometric. For example, every element of
order $5$ is conjugate to a tranformation of the form
$z\mapsto (exp(2k\pi i/q))z$, for $k=1,\ldots,4$, but only those with
$k=1,4$ are geometric.

An element $A$ of a group $G$ of \M transformations is called {\bf primitive}
if the only solutions of the equation $B^n=A$, with $B\in G$ and
$n\in\Bbb Z$, are given by $B=A^{\pm 1}$.

\sn {\bf 1.2.} Let $G$ be a non-elementary finitely generated
Kleinian group. An isomorphism
$\theta:G\rightarrow\theta(G)\subset{\rm PSL}(2,\Bbb C)$ is called
{\bf geometric} if
there exists a quasiconformal homeomorphism of the Riemann sphere $w$,
such that $\theta(g)=w\circ g\circ w^{-1}$ for all $g\in G$. Two geometric
isomorphisms
$\theta_1,~\theta_2$, are equivalent if there exists a \M
transformation $A$, such that $\theta_1(g)=A\circ\theta_2(g)\circ A^{-1}$,
for all $g\in G$.
The set of equivalence classes of geometric isomorphisms of $G$ is
the {\bf \T} (or {\bf deformation}) space of $G$, $T(G)$. It is a well known
fact
that $T(G)$ is a complex manifold (\cite{bers:spaces}, \cite{kra:spaces},
\cite{mas:self}).

\sn {\bf 1.3.} A {\bf signature} is a set of numbers
$\sigma=(p,n;\nu_1,\ldots,\nu_n)$ such that
$p,~n\in {\Bbb Z}^+$, $\nu_j\in {\Bbb Z}^+\cup\{\infty\}$,
$\nu_j\geq 2$. We will say the signature is {\bf hyperbolic}, {\bf parabolic}
or {\bf elliptic} if the quantity $2p-2+n-\sum_1^n(1/\nu_j)$ is positive, zero
or negative, respectively. The pair $(p,n)$ is called the {\bf type} of the
signature.

A ($2$-){\bf orbifold} $\cal S$ of signature $\sigma$, is a
compact surface of genus $p$ from which finitely many points have been
removed (as many as $\infty'$s are in $\sigma$) and such that $\cal S$
has a covering which is locally $\nu_j$-to-1 over certain points, where the
$\nu_j'$s are the finite values appearing in $\sigma$. The $\nu_j'$s
are called {\bf ramification values} of the signature or of the orbifold.
A {\bf maximal partition} $\cal C$ on an orbifold $\cal S$ with
hyperbolic
signature $\sigma$ is a set of $3p-3+n$ simple unoriented disjoint curves on
${\cal S}'={\cal S}-\{x_j;~\nu_j<\infty\}$, such that no two curves of
$\cal C$ are
freely homotopically equivalent on $S'$,
and no curve of $\cal C$ is freely
homotopically equivalent to a loop around a point or a
puncture of ${\cal S}'$.
To avoid trivial cases, when we talk about maximal partitions
we will assume that the signature of the orbifold
satisfies $3p-3+n>0$.

\sn {\bf 1.4.} The following result of B. Maskit provides us with a
uniformization of orbifolds by Kleinian groups that are better for
computational purposes than Fuchsian groups.
\setcounter{thm}{0}
\begin{thm}[Maskit \cite{mas:bound}, \cite{mas:class}]
Given an orbifold $\cal S$ with hyperbolic signature
$\sigma$ and maximal partition $\cal C$,
there exists a (unique up to conjugation in {\rm PSL}($2,\Bbb C$))
geometrically finite Kleinian group $\Gamma$,
called a {\bf terminal regular b-group}, such that:

1.- $\Delta/\Gamma$ is conformally equivalent to $\cal S$;

2.- for each element $a_j$ of the partition $\cal C$, there is a curve
$\tilde{a}_j$ in $\Delta$, precisely invariant under a cyclic subgroup
$<A_j>$ of $G$, generated by an accidental parabolic transformation $A_j$,
and such that
$\Delta \supseteq \tilde{a}_j\stackrel{\pi}{\rightarrow}
\pi(\tilde{a_j})=a_j\subseteq
{\cal S},$ where $\pi:\Delta\rightarrow\cal S$ is the natural projection.

3.- $(\Omega (\Gamma)-\Delta)/\Gamma$ is the union of the orbifolds of type
$(0,3)$ obtained by squeezing each curve of $\cal C$ to a puncture
and discarding all orbifolds of signature $(0,3;2,2,\infty)$
that appear.
\end{thm}

{}From a \T theory point of view, the only interesting surface
uniformized by a terminal regular b-group $\Gamma$ is the one corresponding to
the invariant component, since the deformation spaces of
orbifolds of type $(0,3)$ are points. The space $T(\Gamma)$
is then equivalent to the deformation space of the orbifold $T(S)$.
When convenient, we will denote
$T(\Gamma)$ by $T(p,n;\nu_1,\ldots,\nu_n)$. Its complex dimension is
$3p-3+n$.

\sn {\bf 1.5.} The uniformization theorem of Maskit allows us to embed
$T(\Gamma)$ into
a product of one dimensional \T spaces as follows. Let $T_j$ be the connected
component of ${\cal S}-\{a_k;~a_k\in{\cal C},~k\neq j\}$ containing the curve
$a_j$.
Let $D_j$ be a connected component of $\pi^{-1}(T_j)$, and let
$\Gamma_j$ be its stabilizer in $\Gamma$;
that is, $\Gamma_j=\{\gamma\in \Gamma;~\gamma(D_j)=D_j\}$. The
$\Gamma_j$'s are terminal regular b-groups of type $(0,4)$ or $(1,1)$,
and therefore $T(\Gamma_j)$ is a one-dimensional manifold
(\cite{kra:variational}).
It is clear that any geometric isomorphism
of $\Gamma$ induces an geometric isomorphism of $\Gamma_j$ by restriction.
\begin{thm}[Maskit \cite{mas:moduli},
Kra \cite{kra:variational}]
The mapping defined by restriction,
$T(\Gamma)\hookrightarrow \prod_{j=1}^{3p-3+n} T(\Gamma_j)$,
is holomorphic and injective with open image.
\end{thm}
\noindent We will call this mapping the {\bf Maskit embedding}
of the group $\Gamma$.

\sn {\bf 1.6.} Throughout this paper, for a signature
$(p,n;\nu_1,\ldots,\nu_n)$, $q_j$ and $p_j$ will denote
$\cos (\pi/\nu_j)$ and $\sin (\pi/\nu_j)$, respectively, $j=1,\ldots,n$.

\section{Triangle Groups}
\noindent {\bf 2.1.} A triangle groups is a Kleinian group $\Gamma$
generated by two
M\"{o}bius transformations, $A$ and $B$, such that $|A|=\nu_1$,
$|B|=\nu_2$ and $|AB|=\nu_3$. Here $|T|$ denotes the order of the
transformation $T$, with parabolic elements considered as elements of
order equal to $\infty$. A triangle group $\Gamma$ with hyperbolic
signature can be constructed, for example, by taking
a triangle on $\Bbb H$ with angles $\pi/\nu_j$, $j=1,2,3$,
and considering the group of transformations $\Gamma^*$
generated by reflections on the sides of the triangle.
Then $\Gamma$ is the index $2$ subgroup consisting on the orientation
preserving transformations. A similar construction can be carried
out for the case of parabolic groups. Our goal is to correctly choose
the position of the vertices of such triangles.

{\large\bf\noindent Hyperbolic groups.}\\ \smallskip
\sn {\bf 2.2.} Let $(a,b,c)$ be
three distinct points on \C. Let $\Lambda$ be the
circle determined by these points. Let $\Delta=\{z\in\hat{\Bbb
C};~{\rm Im}(cr(a,b,c,z))>0\}$, where $cr$ denotes the cross ratio of four
different points of the Riemann sphere, chosen such that
$cr(\infty,0,1z)=z$ (remember that there are six possible definitions
of cross ratio). Let $L$ and $L'$ be the circles orthogonal to
$\Lambda$ and passing through $\{a,b\}$ and $\{a,c\}$, respectively.
\begin{defn} Let $z_1$ and $z_2$ be two distinct points in
$L\cap\overline{\Delta}.$ We will say that
they are WELL ORDERED, with respect to ($a,b,c$),
if one of the following
set of conditions is satisfied (they are not mutually exclusive):

{\rm 1}) $z_1=a$;

{\rm 2}) $z_2=b$;

{\rm 3}) $z_1\neq a$, $z_2\neq b$ and $cr(a,z_1,z_2,b)$ is real and
strictly bigger than $1$.
\end{defn}
For example, if $a=\infty$, $b=0$ and $c=1$, we have that $L$ is the
imaginary axis, and $\Delta$ is the upper half plane. If $z_1=\lambda
i$ and $z_2=\nu i$, then these two points are well ordered with
respect to $(\infty,0,1)$ if and only if $\lambda>\nu>0$.\\
Given this definition, we can state the concept of \lq good\rq\ generators.
\begin{defn} Let $(a,b,c)$ be three distinct points of \C,
and let $\Lambda,~\Delta,~L$ and $L'$ be as previously defined.
Suppose that $\Gamma$ is a triangle group with hyperbolic
signature $(0,3;\nu_1,\nu_2,\nu_3)$ and whose limit set is $\Lambda$.
Let $A$ and $B$ be elements of $\Gamma$.
We will say that $(A,B)$ are CANONICAL GENERATORS of $\Gamma$
for the parameters
$(a,b,c)$ if they generate the group $\Gamma$ and the following conditions are
satisfied:

{\rm 1}) $|A|=\nu_1$, $|B|=\nu_2$, $|AB|=\nu_3$,

{\rm 2}) $A$ and $B$ have their fixed points on $L$, and $AB$ on $L'$,

{\rm 3}) if $z_1$ and $z_2$ are the fixed points of $A$ and $B$ on
$L\cap\overline\Delta$, then they are well ordered with respect
to ($a,b,c$),

{\rm 4}) $A$ and $B$ are primitive elements, and geometric whenever elliptic.
\end{defn}

Our main result about existence and uniqueness of canonical
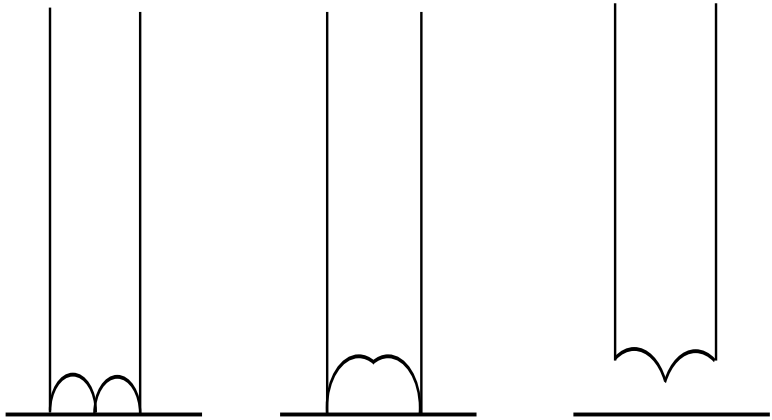
generators for hyperbolic triangle groups is the following:
\begin{thm}Given three different points $(a,b,c)$ in the Riemann sphere,
and a hyperbolic signature $\sigma=(0,3;\nu_1,\nu_2,\nu_3)$,
there exists a unqie pair of \M transformations ($A,B$), such that
they are canonical generators of a triangle group with
signature $\sigma$ and for the given parameters.\\
In the case $(a=\infty,b=0,c=1)$, these generators are given by:

{\rm 1}) Signature $(0,3;\infty,\nu_1,\nu_2)$, $\nu_i\leq\infty$:
$$A=\left[ \begin{array}{cc} -1&-2\\0&-1\end{array}\right],~
B=\left[ \begin{array}{cc} -q_1&b\\q_1+q_2&-q_1\end{array}\right],~
b=\frac{q_1^2-1}{q_1+q_2}$$.

{\rm 2}) Signature $(0,3;\nu_1,\nu_2,\nu_3)$, where all the $\nu_i$
are finite:
$$A=\left[ \begin{array}{cc} -q_1&-kp_1\\k^{-1}p_1&-q_1.\end{array}\right],~
B=\left[ \begin{array}{cc} -q_2&-hp_2\\h^{-1}p_2&-q_2\end{array}\right].$$
Here
$k=\frac{q_2+q_1q_3+q_1l}{p_1l},~h=\frac{kp_1p_2}{q_1q_2+q_3+l}$,
and
$l=\sqrt{q_1^2+q_2^2+q_3^2+2q_1q_2q_3-1}$
with the square root chosen to be positive.\\
For any other set of parameters, $(a,b,c)$, the generators are
given by conjugating the
above transformations by the unique \M transformation $T$ that maps
$a,b,c$ to $\infty,0,1$ respectively.\end{thm}

Observe that the transformations in the second case of the above
theorem converge to those of the first case when $\nu_1$ goes to $\infty$.

\sn {\bf 2.3.} The following two technical results are needed to
prove theorem $3$.
\begin{prop}
Let $\Gamma$ be a triangle group with signature $(0,3;\mu_1,\mu_2,\mu_3)$,
Let $l^2=q_1^2+q_2^2+q_3^2+2q_1q_2q_3-1.$
Then $l^2$ is positive, zero or negative if and only if the signature
is hyperbolic, parabolic or elliptic, respectively.
\end{prop}
\begin{pf} The elliptic and parabolic cases can be checked by computing the
values of $l^2$. In the hyperbolic case, first observe that the
expression of $l^2$ is symmetric on $\nu_j'$s, so we can assume,
without loss of generality, that $\nu_1\leq\nu_2\leq\nu_3$. It is also
clear that $l^2$ is increasing with $\nu_j$, so we need to compute
its values only for the cases of small signatures. More precisely,
it suffices to consider the cases $(0,3;2,3,7)$ and
$(0,3;3,3,4)$. For the first of these signature we have
$l^2=\cos^2(\pi/7)-\frac{3}{4}>\cos^2(\pi/6)-\frac{3}{4}=0$, while
for the second signature we get $l^2=1+\frac{3\sqrt{2}}{4}>0$.\end{pf}
\begin{prop}Let $A$ and $B$ be canonical generators for the group
$\Gamma(\nu_1,\nu_2,\nu_3;a,b,c)$.
Assume that $\tilde{A}$ and $\tilde{B}$ are liftings of $A$ and $B$,
respectively, to SL($2,{\Bbb C}$), both having negative trace. Then
the product $\tilde{A}\tilde{B}$ also has negative trace.
\end{prop}
\begin{pf}
We start with the observation that all the ramification values should
be bigger than $2$, since involutions have matrix representatives
with zero trace, and then the proposition would not make sense. Since
the trace of a matrix is invariant under conjugation, we are free to
choose $(\infty,0,1)$ as parameters for the group; this will simplify
our computations.

Let us first look at the case of
$\Gamma(\infty,\nu_2,\nu_3;\infty,0,1)$. Assume that $\tilde{A}$ and
$\tilde{B}$ have negative trace but $\tilde{A}\tilde{B}$ has positive
trace. We have the following expressions for the liftings of $A$ and $B$:
$$\tilde{A}=\left( \begin{array}{cc} -1&-\alpha\\0&-1\end{array}\right),~
\tilde{B}=\left( \begin{array}{cc} a&b\\c&-a-2q_2\end{array}\right).$$
{}From the definition of canonical generators for the above parameters
we get the following conditions:
$$\left\{\begin{array}{lcl}
{\rm trace}(\tilde{A}\tilde{B})=-2q_1\leq 0
&\Leftrightarrow&2q_2-\alpha c=2q_3\\
\mbox{\rm Re(fixed points of B)}=0
&\Leftrightarrow&2a+2q_2=0\\
\mbox{\rm Re(fixed points of AB)}=1&
\Leftrightarrow&-\alpha c-2q_2-2a=-2c
\end{array}\right .$$
Solving these equations we get
$$\tilde{A}=\left( \begin{array}{cc} -1&-2\\0&-1\end{array}\right),~
\tilde{B}=\left( \begin{array}{cc} -q_1&b\\q_1-q_2&-q_1\end{array}\right),
b=\frac{q_1^2-1}{q_1-q_2}.$$
The Schimizu-Leutbecher Lemma \cite[pg. 18]{mas:kg}
implies that $|q_2-q_3|\geq
1/2$ or $q_2-q_3=0$. The second situation cannot happen as it would
imply that the group is elementary.
In the first case, the only
possible solutions are $q_2=1$, $q_3=1/2$ or vice versa. This implies
that the possible signatures of the group are $(0,3;\infty,\infty,3)$
or $(0,3;\infty,3,\infty)$. But in both cases we obtain an element of
order $2$, namely $ABA$ and $AB^{-1}$ respectively, which is not
possible because of the signatures.

The case of the group $\Gamma(\nu_1,\infty,\nu_3;\infty,0,1)$, with
$\nu_1<\infty$, can be reduced to the previous situation as follows.
Let $x\in {\Bbb R}-\{0\}$ be the end point of the
geodesic in $\Bbb H$ joining
$0$ and the fixed point of $AB$. Then
$B$ and $(AB)^{-1}$ are canonical generators for
$\Gamma(\infty,\nu_3,\nu_1;0,x,\infty)$. Let $T$ be the \M
transformation that takes $0,x,\infty$ to $\infty,0,1$ respectively.
Then
$T\Gamma(\infty,\nu_3,\nu_1;0,x,\infty)T^{-1}=
\Gamma(\infty,\nu_3,\nu_1;\infty,0,1)$,
and we are in the situation already discussed.

Consider now the case of $\Gamma(\nu_1,\nu_2,\nu_3;\infty,0,1)$,
where all the ramification values are finite. Assume again that
$\tilde{A}\tilde{B}$ has negative trace. Then an easy computation
shows that the matrices representatives of $A$ and $B$ are given by
the following expressions:
$$\tilde{A}=\left( \begin{array}{cc} -q_1&-mp_1\\m^{-1}p_1&-q_1
\end{array}\right),~
\tilde{B}=\left( \begin{array}{cc} -q_2&-np_2\\n^{-1}p_2&-q_2
\end{array}\right),$$
where $m=\frac{q_2-q_1q_3+q_1r}{p_1r}$,
$n=\frac{mp_1p_2}{q_1q_2-q_3+r}$, and $r$ is the positive square root
of \\$q_1^2+q_2^2+q_3^2-2q_1q_2q_3-1$.

If $m>0$, then consider the triangle of figure $1$,
where the vertices $v_1$, $v_2$ and $v_3$ are fixed by $A$, $B$ and
$AB$, respectively, and the angle at $v_j$ is $\pi/\nu_j$.
It is not hard to see that
$A(\infty)<0$ and $AB(\infty)<1$, so the action of the
transformations $A$ and $AB$ are as indicated in the figure. Reflect
the triangle on the geodesic joining $v_1$ and $v_3$ to get a similar
triangle with vertices $v_1$, $v_3$ and $v_4$. Since $A$ and $B$ are
isometries in $\Bbb H$, a look at the figure shows that
$ABA$ fixes $v_2$. This implies that $ABA=B^n$ for some integer $n$.
If $n=0$, we have that the group is elementary, which is not
possible. If $n\neq 0$, we use the fact that $AB=(B^{-1}A^{-1})^{\nu_3-1}$,
to get
$B^{n}=ABA=(B^{-1}A^{-1})^{\nu_3-1}A=(B^{-1}A^{-1})^{\nu_3-2}B^{-1}$.
Therefore $B^{n+1}=(AB)^{\nu_3-2}$. But this equality can be
satisfied only if $\nu_3=2$, contrary to our assumptions.

The case of $m<0$ is solved as above
by considering the element $A^{-1}B^{-1}A^{-1}$ instead of $ABA$.
\end{pf}

\sn {\bf 2.4.} We are now in a position to prove theorem $3$.\\
{\it Proof of theorem $3$}.
\underline{First case}: the signature is $(0,3;\infty,\nu_1,\nu_2)$.
By the same trick used in the proof of proposition $2$,
we can reduce all the signatures with punctures to this case.
The element $A$ is a parabolic transformation
fixing $\infty$, so it is of the form $A(z)=z+\alpha$,
with $\alpha\in{\Bbb C}-\{0\}.$
Consider now the element $B(z)=\frac{az+b}{cz+d}$. Choose liftings of
$A$ and $B$ to SL($2,\Bbb C$) with negative traces.
We then have the following equations:
$$\left\{\begin{array}{lcl}
{\rm trace}(\tilde{B})=-2q_2\leq 0&\Leftrightarrow&a+d=-2q_2,\\
\mbox{\rm Re(fixed points of B)}=0,&
\Leftrightarrow&a-d=0,\\
{\rm trace}(\tilde{A}\tilde{B})=-2q_3\leq 0&\Leftrightarrow&
a-\alpha c-d=2q_3,\\
\mbox{\rm Re(fixed points of AB)}=1&\Leftrightarrow&
a-\alpha c c+d=-2g,\\
{\rm det}\tilde{B}=1&\Leftrightarrow&ad-bc=1.
\end{array}\right .$$
Solving these equations we get the matrices of theorem $1$.

\underline{Second case}: the signature is $(0,3;\nu_1,\nu_2,\nu_3)$,
where all the ramification values are finite.
The canonical generators will be of the form
$A(z)=\frac{\alpha z+\beta}{\gamma z+\delta}$ and
$B(z)=\frac{az+b}{cz+d}$. Choosing matrices representatives of these
transformations with negative traces, we get the equations:
$$\left\{\begin{array}{lcl}
{\rm trace}(\tilde{A})=-2q_1\leq 0&\Leftrightarrow&
\alpha+\delta=-2q_1,\\
\mbox{\rm Re(fixed points of A)}=0&\Leftrightarrow&
\alpha-\delta=0,\\
{\rm trace}(\tilde{B})=-2q_2\leq 0&\Leftrightarrow&
a+d=-2q_2,\\
\mbox{\rm Re(fixed points of B)}=0&\Leftrightarrow&
a-d=0,\\
{\rm trace}(\tilde{A}\tilde{B})=-2q_3\leq 0&\Leftrightarrow&
\alpha a+\beta c+\gamma b+\delta d=-2q_3,\\
\mbox{\rm Re(fixed points of AB)}=1&\Leftrightarrow&
\alpha a+\beta c-\gamma b-\delta d=2(\gamma a+\delta c).
\end{array}\right .$$
Using the fact that the matrices involved in these equations have
determinant equal to $1$, it is not hard to see that the only
solution is the one given in theorem $3$.~$\hfill\Box$

\sn {\bf 2.5.} The following technical result will be neede in \S $2.7$.
\begin{prop}
If ($A,B$) and ($A,D$) are two pairs of canonical generators for a
hyperbolic triangle group with signature $(0,3;\infty,\nu_1,\nu_2)$
then there exists an integer number, $n$, such that $D=A^{n/2}BA^{-n/2}.$
\end{prop}
In general a M\"{o}bius transformation has several square roots, but
in the case of parabolic elements such a situation does not happen.
So we have that
the transformation $A^{n/2}$ is well defined for any integer $n$.

\begin{pf} By conjugation we may assume that
$A$ and $B$ are canonical generators for the parameters $(\infty,0,1)$,
whose expressions are given in theorem $1$.
Let $\tilde{A}$ and $\tilde{D}$ be liftings of $A$ and $D$ to
PSL($2,\Bbb C$) respectively, with negative trace. By proposition $2$
we have that $\tilde{A}\tilde{D}$ has also negative trace.
Computations show that, under these conditions, $A(z)=z+2$ and
$D(z)={\load{\large}{\displaystyle}
\frac{\alpha z+\beta}{(q_1+q_2)z-2q_1-\alpha}}$, $z\in\Bbb C$,
$\alpha,\beta\in\Bbb R$.
Let $T(z)=z+h$, where $h={\load{\large}{\displaystyle}\frac{\alpha+q_1}
{q_1+q_2}}$ (this is the real part
of the fixed points of $AD$).Then $TBT^{-1}=D$.
So $T$ belongs to the normalizer of $\Gamma$ in \MM , and
induces an automorphism of the quotient surface that fixes one puncture (since
$TAT^{-1}=A$). This means that either $T\in\Gamma$ or $T^2\in\Gamma$, giving us
$T=A^n$ or $T=A^{n/2}$ as desired.
\end{pf}\smallskip

{\large\bf\noindent Parabolic groups.}\\ \smallskip
\noindent {\bf 2.6.} The only case of parabolic triangle groups
that we need in this paper
is the one of groups with signature $(0,3;\infty,2,2)$. For a
treatment of the general case, as well as the elliptic signatures,
see \cite{ares:horoc2}.

\begin{defn}
Let $\Gamma$ be a triangle group with signature $(0,3;\infty,2,2)$.
Let $A$ and $B$ be two
generators of the group. We will say that they are CANONICAL for the
parameters $(a,b,c)$ if the following conditions are satisfied:

{\rm 1}) $|A|=\infty$, $|B|=2$ and $|AB|=2$,

{rm 2}) $A(a)=a$, $B(b)=b$ and $AB(c)=c$,

{\rm 3}) $A$ and $B$ are primitive.
\end{defn}
\begin{thm}Given three different points $(a,b,c)$ in the Riemann sphere,
there exists a unique pair of \M transformations ($A,B$), such that
they are the canonical generators of a triangle group with
signature $(0,3;\infty,2,2)$ and for the given parameters. We will
denote the triangle group with the pair of canonical generators by
$\Gamma(\infty,2,2;a,b,c)$.\end{thm}

\begin{pf} Taking parameters $(\infty,0,1)$, the proof is reduced to
a simple computation.\end{pf}

Proposition $4$ also holds for this type of groups. Since the proof is the
same,
we will not reproduce it again.\smallskip

{\large\bf\noindent The geometry of the quotient orbifolds.}\\ \smallskip
\noindent {\bf 2.7.} Our next goal is to produce coordinates on the orbifolds
uniformized by the groups studied above. We will use
these coordinates to explore the relation between parameters on \T
spaces and the construction of Riemann surfaces (see theorems $6$,
and $9$).

Let $\Gamma=\Gamma(\infty,\nu_2,\nu_3;a,b,c)$ be a group with
hyperbolic signature. Let $\Delta$ be defined as in \S 2.1. Put on
$\Delta$ the metric of constant curvature $-1$. Since $\Gamma$ acts
on $\Delta$ by isometries, we have a metric on $S\cong\Delta/\Gamma$
obtained by projection. A geodesic on $S$ is just a curve that lifts
to a geodesic on $\Delta$.
\begin{prop}
Let $\cal S$ be a hyperbolic orbifold with signature
$(0,3;\infty,\nu_1,\nu_2)$. Let $P\not\in\cal S$ be the
puncture
corresponding to the first $\infty$ in the signature, and let $Q$ the
puncture or branched point corresponding to $\nu_1$.
Then there exists a unique
simple geodesic $c:I\longrightarrow\cal S$, joining $P$ and
$Q$ such that, if $c$ is parametrized by arc length $s$, then:

{\rm 1}) if $\nu_1=\infty$, $I=\Bbb R$,
$\lim_{s\rightarrow +\infty}c(s)=P$, and $\lim_{s\rightarrow -\infty}c(s)=Q.$

{\rm 2}) if $\nu_1<\infty$, $I=[0,+\infty)$,
$\lim_{s\rightarrow +\infty}c(s)=P$, and $c(0)=Q$.
\end{prop}
\begin{pf} The existence part is easy. These orbifolds have no moduli, so we
can assume that the covering space is $\Bbb H$,
and the covering group $\Gamma$
has parameters $(\infty,0,1)$. As a fundamental domain for our
group we can choose (depending on the signature of the group)
one of those in figure $2$ below.
The projection of the part of the imaginary axis that lies in
the boundary of that fundamental domain gives a geodesic on
$\cal S$ that satisfies the conditions of the statement of the
proposition.

For the uniqueness part, let us assume that there is another
geodesic, satisfying the properties stated in the proposition.
We lift it to $\Bbb H$ and we can assume that the lifting is
a half-vertical line. We want to
prove that this second vertical line is simply a translate of the imaginary
axis under a power of $A$, and therefore the projection of the two
lines will be the same geodesic in the orbifold.

The end point of our line, say $x_0$, which corresponds to $Q$, has
to be fixed by a transformation $B_1$. Now, if we remove on the
orbifold the point corresponding
to $\nu_2$ (if $\nu_2=\infty$, then we do not have to remove anything,
since punctures are not in the orbifold), we are in a situation
like the torsion free case, and we get that
$A$ and $B_1$ generate the group $\Gamma$ (see \cite{kra:horoc}).
Therefore $A$ and $B_1$ will be
canonical generators for some parameters. By proposition $4$,
we have that there is an integer $n\in{\Bbb Z}$,
such that $A^{n/2}BA^{-n/2}=B_1$.
Our proof will be complete if we show that $n$ is even.

If $n$ is odd, then $B_1=A^{1/2}BA^{-1/2}$ is conjugate to $B$
in the group $\Gamma$, since both transformations, $B$ and $B_1$,
correspond to the same point $Q$ in $S$.
This implies that the element $A^{1/2}$ belongs to the
normalizer of $\Gamma$ in PSL($2,{\Bbb R})$ and therefore it
induces an automorphism of the quotient orbifold
that fixes at least one puncture (the one represented by $A$).
Since $A^{1/2}$ does not belong to $\Gamma$, the induced automorphism
is not the identity, and so it has to interchange the
other two ramification points. This implies that $B$ and $AB$ are
tranformations of the same order. It is easy to see that, in such case,
$A^{1/2}BA^{-1/2}=(AB)^{-1}$. Then we would have that
$B$ and $(AB)^{-1}$ are conjugate in the group $\Gamma$,
which is not true since they correspond to different
branch points. Therefore $x_0$ is an even integer and the geodesics
are the same.\end{pf}

\noindent {\bf 2.8.} We can use the geodesic of the above proposition
to construct coordinates
on the orbifold $\cal S$. In this paper we will use only coordinates
around punctures, and the main result in that line is in \cite{kra:horoc}.
His proof can be applied to our situation
since the argument is local. We copy the result here for the
convenience of the reader.
For a more general situation see \cite{ares:horoc2}
\begin{prop}
Let $\Gamma(\infty,\nu_1,\nu_2;a,b,c)$ be a hyperbolic triangle group.
Let $\cal S$ be an orbifold
uniformized by this group, and suppose that $\cal S$ has the metric
of constant curvature $-1$ that comes from its universal (branched)
covering space. Let $P^1\not\in\cal S$  the puncture corresponding to the
$\infty$ in the signature, and let $P^2$ be the ramification
point corresponding to $\nu_1$. Let $c$ be the geodesic on $\cal S$
given in proposition $4$. Then there exists a biholomorphism $z$, defined in
a neighborhood $N$ of $P^1$, such that $z(P^1)=0$ and $z$ maps
isometrically the portion
of $c$ inside $N$ into the positive real axis, with the metric of the
punctured disc. These properties define the germ of $z$ uniquely.\end{prop}

\sn {\bf 2.9.} We have similar results for the parabolic case.
\begin{prop}
Let $\cal S$ be an orbifold of signature  $(0,3;\infty,2,2)$. Let $P\not\in
\cal S$ be the puncture, and let $Q$ be one of the branched points.
Then there exists a unique geodesic $c:[0,\infty)\longrightarrow\cal S$,
such that $c(0)=Q$, $lim_{s\rightarrow\infty}c(s)=P$, for $s$ the arc
length parametrization, and $c$ realizes the distance between any two
points in it.
\end{prop}
\begin{pf} We first note that there is no loss of generality in
taking the group with parameters
$(\infty,0,1)$, and in assuming that $Q$ lifts to $0$.
By our definition of geodesics, any straight line joining
zero and infinity will project onto a geodesic of the orbifold.
It is easy to see that the imaginary axis $i\Bbb R$ satisfies all the required
properties.\\
To prove uniqueness,
suppose first that the geodesic on $S$ lifts to another vertical line,
say $L$, whose point of intersection with the real axis is $x_0$. Then,
since $x_0$ and $0$ project to the same point on $S$, we must have
that $x_0=2n$, for some integer $n$. This implies that $L$ is the
image of the imaginary axis under the mapping $A^n$, and therefore
$i\Bbb R$ and $L$ give the same geodesic on $S$.\\
To study the cases of non vertical lines,
identify the complex plane with ${\Bbb R}^2$, with coordinates
$(x,y)$. Then we can write the lifting of our geodesic as $L'=\{y=mx\}$, with
$m$ a real number. The slope $m$ cannot be zero because the real
axis projects onto a
line the joins that two branch points, but it stays at away from the
puncture. So we have that $m\neq 0$. Consider the points $2i$ on the
imaginary axis and $2+2im$ on $L'$. Both points project onto the same
point in the orbifold, since they are equivalent under the
transformation $A$.
The distance from $0$ to $2i$ along $i\Bbb R$ is $2$, while
the distance from $0$ to $2+2im$ along $L'$ is $2\sqrt{1+m^2}$.
Therefore $L'$ does not satisfies the properties of the proposition.\end{pf}

\sn {\bf 2.10.} We can construct local coordinates on the
quotient orbifold as in the hyperbolic case,
but we do not have a defining property as the one in
proposition $8$, due to the curvature constrain.
The expression $f_{12}(z)=e^{\pi i\rho^{-1}(z)}$ gives
a local coordinate on ${\cal S}={\Bbb C}/\Gamma(\infty,2,2;\infty,0,1)$
around the puncture. Here $\rho:{\Bbb C}\rightarrow {\cal S}$ is the
natural projection mappping. We will say that the germ
of holomorphic functions defined by $f_{12}$ is a {\bf horocyclic coordinate}.

\section{Coordinates for the Teichm\"{u}ller spaces of b-groups with torsion.}
\noindent {\bf 3.1.} This section is devoted to computing coordinates
for the \T spaces of orbifolds with hyperbolic signature.
As we said in the introduction, to parametrize the \T spaces of
orbifolds, it suffices to consider in detail the one-dimensional
cases. These correspond to the orbifolds of type $(0,4)$ and $(1,1)$.
Then, the Maskit Embedding Theorem gives coordinates in the general
deformation space.

Before starting our computations we need to have a convention
about orientation of curves on orbifolds. Assume that $a$ is a simple
loop contractible to a puncture on an orbifold $S$. Then we
orient $a$ by requiring that the puncture lies to the left of the curve.

Similarly, if $A$ is a parabolic transformation fixing $z_0\in\hat{\Bbb
C}$, and $L$ is a horocircle through $z_0$, we orient $L$ by
choosing a point $z\neq z_0$ in $L$, and requiring that $z$, $A(z)$
and $A^2(z)$ follow each other in the positive orientation.
Observe that a horocircle passing through $\infty$ can be undertood
as a circle on the Riemann sphere, and therefore it makes sense to
talk about its positive orientation.

\sn {\bf 3.2.} Let us start with the case of an orbifold $S$ of signature
$(0,4;\nu_1,\ldots,\nu_4)$. A maximal partition on $S$ consists of a
curve $a$, that divides $S$ into two subsets $S_1$ and $S_2$, each of
them with two ramification points. Without loss of generality, we can
assume that the points with ramification values $\nu_1$ and $\nu_2$
are in $S_1$, and this set lies to the right of $a$. Cut $S$ along
$a$ and glue to each resulting boundary curve a punctured disc. In
this way we complete $S_j$ to an orbifold of type $(0,3)$. Denote
these new orbifolds by the same letter, and assume for the moment
that both of them, $S_1$ and $S_2$, are of hyperbolic type.
We have that $S_1$ is uniformized by a triangle group, that can be taken to
be $\Gamma_1=\Gamma(\infty,\nu_1,\nu_2;\infty,0,1)$. Since we are
interested on parameters for \T spaces, we are free to conjugate by a
\M transformation, which explains our choice of $\Gamma_1$.
Its canonical generators are given in theorem $3$.
Since $S_1$ lies to the right of $a$, our orientation requirements
imply that we are considering the action of $\Gamma_1$ on the upper
half plane.

The orbifold $S_2$ is uniformized by a triangle group
$\Gamma_2=\Gamma(\infty,\nu_3,\nu_4;d,e,f)$.
We have that the transformation that corresponds to the curve $a$ is
the canonical generator $A(z)=z+2$. Since $S_1$ and $S_2$ come both
from the same orbifold, and we want to glue them to obtain $S$,
we must have that the element that corresponds to $a$ in $\Gamma_2$
must be either $A$ or $A^{-1}$. But $S_2$ lies to the left of the
partition curve, so we get that one of the canonical generators of
$\Gamma_2$ is $A^{-1}$. This implies that $\Gamma_2$ is conjugate to
$\Gamma(\infty,\nu_3,\nu_4;\infty,0,-1)$ by a transformation $T$, such
that $TA^{-1}T^{-1}=A^{-1}$. Therefore $T(z)=z+\alpha$,
$\alpha\in{\Bbb C}-\{0\}$. The generators of $\Gamma_2$ are then
given by $A^{-1}$ and
$$B_\alpha^{-1}=\left[\begin{array}{cc}-q_3-\alpha (q_3+q_4)&
-b^*-\alpha^2 (q_3+q_4)\\
-(q_3+q_4)&-q_3+\alpha (q_3+q_4)\end{array}\right],~
b^*=\frac{1-q_3^2}{q_3+q_4}.$$

We have ${\rm Im}(\alpha)>0$, since $S_1$ is uniformized by
$\Gamma_1$ in the upper half plane. Suppose
now that this imaginary part is big enough, say bigger than $2$. Then
the curve $\{{\rm Im}(z)={\rm Im}(\alpha)/2\}$ is invariant under $A$
and lies in the intersection of the discontinuity sets of $\Gamma_1$
and $\Gamma_2$. So we can use this curve to apply the
First Klein-Maskit Combination Theorem \cite[pg. 149]{mas:kg}. In this way we
obtain that the group
$\Gamma_\alpha=\Gamma_1*_{<A>}\Gamma_2:=<\Gamma_1,\Gamma_2>$ is a
terminal
regular b-group uniformizing an orbifold with the above signature. By the
classical theory of quasiconformal mappings, we have that any orbifold
of the above signature is uniformized by one such group $\Gamma_\alpha$.

The parameter $\alpha$ is then a global coordinate for the space\\
$T(0,4;\nu_1,\ldots,\nu_4)$. It can be expressed in an invariant way
as follows. Suppose
$\Gamma=\Gamma(\infty,\nu_1,\nu_2;a,b,c)*_{<C>}
\Gamma(\infty,\nu_3,\nu_4;d,e,f)$ is a terminal regular b-group uniformizing
an orbifold
with the above signature. Then the point corresponding to $\Gamma$ in
the deformation space is $\alpha=cr(a,b,c,f)$.

\begin{thm}$\alpha$ is a global coordinate, called a {\bf horocyclic
coordinate}, for\\
$T(0,4;\nu_1,\nu_2,\nu_3,\nu_4)$. The following inclusions are satisfied:
$$\{\alpha;~{\rm Im}(\alpha )> y_1\}\subset
T(0,4;\nu_1,\nu_2,\nu_3,\nu_4)\subset
\{\alpha;~{\rm Im}(\alpha )> y_2\},$$
where
$$y_1=\frac{1}{q_1+q_2}+\frac{1}{q_3+q_4},\;
y_2=max(\frac{1}{q_1+q_2},\frac{1}{q_3+q_4}).$$
\end{thm}
\begin{pf} We must show only the inclusions.
The first one follows from the above reasoning about
how to apply the First Combination Theorem.

For the second inclusion, we must use the fact that the lower half
plane  $X$ is precisely invariant under $\Gamma_1$ in $\Gamma_\alpha$;
that is, if $\gamma\in\Gamma_\alpha$ and $\gamma(X)\cap X\neq\emptyset$,
then $\gamma\in\Gamma_1$. The point
$z=\frac{-q_3-\alpha(q_3+q_4)-i}{q_3+q_4}$, has negative
imaginary part; and the imaginary part of $B_\alpha(z)$ is equal to ${\rm
Im}(\alpha)-\frac{1}{q_3+q_4}$. This last number should be positive,
giving one condition on the imaginary part of $\alpha$\\
Similarly, we have that the set $Y=\{z;~{\rm Im}(z)>{\rm Im}(\alpha)\}$
is precisely invariant under $\Gamma_2$ in $\Gamma_\alpha$.
Consider the point $w=\frac{q_1+i}{q_1+q_2}$. Its image under $B$ has
imaginary part equal to $\frac{1}{q_1+q_2}$. So, if ${\rm
Im}(\alpha)<\frac{1}{q_1+q_2}$, we then get that both, $w$ and $B(w)$,
belong to $Y$.\end{pf}

\sn {\bf 3.3.} We can relate the above group theoretical
computations to a more geometric construction by means of the
plumbing contructions. This is a well known technique, and we will
not explain it here in detail. See \cite{kra:horoc} for a careful
treatment of it. We will only say that, given two orbifolds $S_1$
and $S_2$, with punctures $P_1$ and $P_2$ respectively, one can construct
a new orbifold as follows.
Remove neighborhoods $V_j$ of $P_j$, $j=1,2$, and identify the
boundaries of the resulting orbifolds $S_j-V_j$. If such
identification is given by am expression of the form $z_1z_2=t$,
where $z_j$ is a horocyclic coordinate around $P_j$, then we say that
the resulting orbifold has been constructed by {\bf plumbing}
techniques with parameter $t$. One can as well do plumbing
contructions in one single orbifold.

To compute the plumbing parameter in the above construction of groups
of type $(0,4)$, for the group $\Gamma_1$ we take the coordinate
$z(\xi)=e^{\pi i\xi}$, with $\xi$ in the upper half plane.
Similarly, for the
orbifold $S_2$ we take $w(\xi)=e^{\pi i(\alpha-\xi)}$, with ${\rm
Im}(\xi)<{\rm Im}(\alpha)$. Then we get $zw=t=e^{\pi i\alpha}$.
We have proven the following results,
up to the (easy) computation on the bounds for plumbing parameters.
\begin{thm}The orbifold corresponding to the point $\alpha$ in
$T(0,4;\nu_1,\nu_2,\nu_3,\nu_4)$ is conformally equivalent
to the orbifold constructed by plumbing with parameter $t=e^{\pi i\alpha}$.
We have that $0<|t|<e^{-\pi y_2},$ with $y_2$ as in theorem $5$.\end{thm}

\sn {\bf 3.4.} The case of one of the orbifolds, say $S_2$, having signature
$(0,3;\infty,2,2)$ is treated in a similar way. We will leave the
computations to the reader, writing only the final results. The group
$\Gamma_1$ is given by $\Gamma(\infty,\nu_1,\nu_2;\infty,0,1)$, while
$\Gamma_2=\Gamma(\infty,2,2;\infty,\alpha,\alpha-1)$. The generators
for $\Gamma_1$ are given in theorem $3$; those of $\Gamma_2$
are
$$A^{-1}=\left[\begin{array}{cc}-1&2\\0,-1\end{array}\right],~
B_\alpha=B^{-1}_\alpha=
\left[\begin{array}{cc}-i&2i\alpha\\0&i\end{array}\right] .$$
\begin{thm}$\alpha$ is a global coordinate for
$T(0,4;\nu_1,\nu_2,2,2)$, and
we the following inclusions are satisfied:
$$\{\alpha;~{\rm Im}(\alpha)>0\}\subset
T(\Gamma(0,4;\nu_1\nu_2,2,2))\subset
\{\alpha;~{\rm Im}(\alpha )> y_1\},$$
where $y_1={\load{\large}{\displaystyle}\frac{1}{q_1+q_2}}$.\end{thm}
We also have a result about plumbing parameters similar to the one in
theorem $6$.

\sn {\bf 3.5.} The other deformation spaces of dimension one correspond
to orbifolds
of signature $(1,1;\nu)$. Let $S$ be an orbifold with this signature,
and let $a$ be a maximal partition on $S$. If we cut $S$ along $a$
and glue punctured discs to the boundary curves, what we get is a
single orbifold $S_1$ with signature $(0,3;\infty,\infty,\nu)$. Let
$\Gamma_1=\Gamma(\infty,\infty,\nu;\infty,0,2)$ be a triangle group
uniformizing $S_1$. To obtain $S$, what we have to do is to
identify the elements $A$ and $B$ corresponding to the two punctures
on $S_1$. Due to the orientation of the curves around punctures, the
correct identification is carried by a \M transformation $C$ such
that $CB^{-1}C^{-1}=A$. The expressions of $A$ and $B$ are again
given in theorem $1$. An easy computation shows that $C$ has the form
$$C=\left[\begin{array}{cc}i\tau&i\sqrt{\frac{2}{1+q}}\\
i\sqrt{\frac{1+q}{2}}&0\end{array}\right].$$
Applying the Second Klein-Maskit Combination Theorem we get that the
group
$\Gamma_\tau=\Gamma_1*_C:=<\Gamma_1,C>$ is a terminal regular b-group of the
desired signature.
\begin{thm}$\tau$ is a global coordinate for $T(1,1;\nu)$, and
we have the following inclusions:
$$\{\tau ;~{\rm Im}(\tau )> 2\}\subset
T(\Gamma(1,1;\nu))\subset\{\tau ;~{\rm Im}(\tau )> 0\}.$$
\end{thm}
\begin{pf} Observe that $C$ maps horocircles at
$0$ (that is, circles passing through zero) to straight lines
(horocircles at $\infty$). In particular we have
$$C(\{z;~|z-ri|=r\})
=\{z;~{\rm Im}(z)=\sqrt{\frac{2}{1+q}}{\rm {\rm Im}}(\tau)-
\frac{1}{r(1+q)}\}.$$
If these two circles are disjoint, the Second Combination Theorem can
be applied. Therefore we want
$$\sqrt{\frac{2}{1+q}}{\rm Im}(\tau)-\frac{1}{r(1+q)}>2r,\;{\rm or}\;
{\rm Im}(\tau)>\sqrt{\frac{1+q}{2}}(\frac{1}{r(1+q)}+2r).$$
The minimum value of the last expression is $2$. This gives the first
inclusion of the theorem. For the other inclusion we just need to use
that the lower half plane is precisely invariant under $\Gamma_1$ in
$\Gamma_\tau$. So for any point $z$ with negative imaginary part we
shold have
$${\rm {\rm Im}}(C(z))=\sqrt{\frac{2}{1+q}}{\rm {\rm Im}}(\tau)-\frac{2}{1+q}
\frac{{\rm {\rm Im}}(z)}{|z|}>0,$$
which gives the desired result.\end{pf}

\sn {\bf 3.6.} As in the $(0,4)$ case, we have a plumbing construction
for these
orbifolds. Consider the coordinates on $S_1$ given by $z=e^{\pi
i\zeta}$ and $w=-exp(\frac{-2\pi i}{(1+q)\zeta})$, near the punctures
determined by $A$ and $B$, respectively. The identification we have
to make in this case is given by
$$z(C(\zeta))w(\zeta)=
z(\sqrt{\frac{2}{1+q}}+\frac{2}{1+q}\frac{1}{\zeta})w(\zeta)=
exp(\frac{2\pi i}{1+q})exp(\sqrt{\frac{2}{1+q}}\pi i\tau)=t.$$
This proves the following result:
\begin{thm}The orbifold corresponding to the point $\tau$ in
$T(1,1;\nu)$ is conformally equivalent
to the orbifold constructed by plumbing with parameter $t$,
with $0<|t|<1.$\end{thm}

\sn {\bf 3.7.} Let $S$ now be a hyperbolic orbifold with a maximal
partition $\cal S$,
uniformized by the terminal regular b-group $\Gamma$. As we noted
in \S $3.1$, we can decompose the group into subgroups with one
dimensional deformation spaces. We have explained in detail how to
parametrize these simpler \T spaces. Together with the Maskit
Embedding, we obtain the following results for the general situation.
\begin{thm}Let $\cal S$ be an orbifold of hyperbolic type with
signature $\sigma=(p,n;\nu_1,\\
\ldots,\nu_n)$
and let $\cal C$ be a maximal partition on $\cal S$, uniformized by the
terminal regular b-group $\Gamma$.
Then there exists
a set of (global) coordinates, $(\alpha_1,\ldots,\alpha_d)$,
called {\bf horocyclic coordinates},
for the \T space $T(\Gamma)\cong T(p,n;\nu_1,\ldots,\nu_n)$,
where $d=3p-3+n$, and a set of complex numbers,
$(y_1^1,\ldots,y_1^d,y_2^1,\ldots,y_2^d)$,
which depends
on the signature $\sigma$ and the partition $\cal C$, such that
$$\{(\alpha_1,\ldots,\alpha_d)\in {\Bbb C}^d;~
{\rm Im}(\alpha_i)> y_1^i\; ,\forall \;1\leq i\leq d\}\subset
T(p,n;\nu_1,\ldots,\nu_n)$$
and
$$ T(p,n;\nu_1,\ldots,\nu_n)\subset
\{(\alpha_1,\ldots,\alpha_d)\in {\Bbb C}^d;~
{\rm Im}(\alpha_i)> y_2^i\; ,\forall \;1\leq i\leq d\}.$$
Moreover, the surface corresponding to the point
$(\alpha_1,\ldots,\alpha_d)$ is conformally equivalent
to a surface constructed by plumbing techniques with
parameters $(t_1,\ldots,t_d)$, obtained
as in theorems $5$ and $8$.
\end{thm}
It is clear that given a point $\alpha$ in
$T(p,n;\nu_1,\ldots,\nu_n)$, one can construct a set of M\"obius
transformations that generate a group $\Gamma$, which corresponds to
$\alpha$. The one dimensional case has been done
explicitly. In the general case, one has simply to iterate the
constructions explained above. For a more detailed description of
this process, see \cite[\S $7.5$]{kra:horoc}, with the necessary
modifications to include finite order transformations.
The above techniques can  be generalized to constructions of Kleinian
groups without parabolic elements as well; see \cite{ares:horoc2}.

\section{The Patterson Isomorphisms in the horocyclic coordinates.}
{\bf 4.1.} One of the most natural questions one may ask
about \T spaces is that under
what circumstances two such spaces are biholomorphic.
A result of Patterson (\cite{pat:dist}, \cite{ek:hol}) states that all
possible isomorphisms between \T space of hyperbolic orbifolds of
different type (with $2p-2+n>0$) are
$T(2,0)\cong T(0,6;2,2,2,2,2,2),$
$T(1,2;\infty,\infty)\cong T(0,5;\infty,2,2,2,2)$ and
$T(1,1;\infty)\cong T(0,4;\infty,2,2,2).$
The existence of these isomorphisms is based on the fact that all surfaces of
genus 2, or of genus 1 with either two or one punctures, have a conformal
involution (hyperelliptic involution); the quotient of the surface by that
involution is a sphere with
six, five or four points, with ramification values as above. Our main
result is as follows.
\begin{thm}The mapping
$$(\tau_1,\tau_2,\tau_3)\mapsto (\frac{\tau_1}{2},1+\tau_2,
1+\frac{\tau_3}{2})$$
gives an isomorphism between $T(0,2)$ and $T(0,6;2,\ldots,2)$ in
the horocyclic coordinates $\tau_j$, $j=1,2,3$ corresponding to the
partition given in figure $3$.
\end{thm}

\sn {\bf 4.2.} We start with a surface of genus 2 with a maximal partition as
shown in figure $3$.

Let $\Gamma$ denote a terminal regular b-group uniformizing the surface and
the partition in the simply connected invariant component $\Delta$.
A presentation for $\Gamma$ can be found in \cite{kra:horoc};
we copy it here for the convenience of the reader.
$\Gamma=\{A_1,~C_1,~A_3,~C_3;~A_1,~A_2=[C_1^{-1},A_1],~A_3~
\mbox{\rm are accidental parabolic,}~
[A_1,C_1^{-1}]\circ [A_3^{-1},C_3^{-1}]=I\}$,
where $[A,B]=ABA^{-1}B^{-1}$. The elements $A_i$ correspond to the
curves $a_i$, while $C_i$ correspond to $c_i$.
These transformations have the following espressions:
$$A_1=\left[\begin{array}{cc} -1&-2\\0& -1\end{array}\right],~
A_2=\left[\begin{array}{cc}1&-2\\-2&-3\end{array}\right],~
A_3=\left[\begin{array}{cc}-1-2\tau_2(1-\tau_2)&-2(1-\tau_2)^2
\\2\tau_2^2&
-1+2\tau_2(1-\tau_2)\end{array}\right],$$
$$C_1=i\left[\begin{array}{cc}\tau_1&1\\1&0\end{array}\right],~
C_3=i\left[\begin{array}{cc}\tau_3\tau_2^2+2(1-\tau_3)\tau_2+\tau_3-2&
-\tau_3\tau_2+(3\tau_3-2)\tau_2-2\tau_3+3\\
\tau_3\tau_2+(2-\tau_3)\tau_2-1&-\tau_3\tau_2^-2(1-\tau_3)\tau_2+2
\end{array}\right].$$
$\tau_1,~\tau_2$ and $\tau_3$ are complex numbers.

The key ingredients in the proof of our theorem is the fact that
the hyperelliptic involution lifts to a M\"obius transformation in
the covering determined by $\Gamma$ and $\Delta$. More precisely, we
have that such lifting is given by the transformation
$A_2^{1/2}={\load{\large}{\displaystyle}\frac{1}{-z+2}}$
(see \cite{kra:horoc}).
\begin{prop}The subgroup $\tilde{\Gamma}$ of PSL(2,$\Bbb C$)
generated by $\Gamma$ and $A_2^{1/2}$
is a terminal regular b-group of signature  $(0,6;2,2,2,2,2,2)$.
\end{prop}
\begin{pf}The facts that $\tilde{\Gamma}$ is Kleinian and
geometrically finite follow
from \cite{mas:kg}, V.E.10, in p. 98 and VI.E.6 in p. 132, respectively.
We also get that $\Omega ({\tilde{\Gamma}})=\Omega (\Gamma)=\Omega$.
Recall that $\Delta$ is the invariant component of $\Gamma$.
If $A_2^{1/2}(\Delta)=U$, where
$U$ is some component of $\tilde{\Gamma}$, then
for all $\gamma\in\Gamma$ we have $A_2^{1/2}\gamma A_2^{-1/2}(U)=U$.
This implies that $A_2^{1/2}(\Delta)=\Delta$, since
$A_2^{1/2}\gamma A_2^{-1/2}\in\Gamma$ and $\Delta$ is the unique
invariant component of this group. The signature of
$\Delta/\tilde{\Gamma}$ is a consequence of the fact that
the hyperelliptic involution
fixes six points on the surface $\Delta/\Gamma$.
The statement about the accidental
parabolic elements of $\tilde{\Gamma}$ is trivial.

To finish the proof we need to show that $(\Omega-\Delta)/\tilde{\Gamma}$ is a
union of orbifolds of type $(0,3)$ (of certain signatures). Let $\Omega_0$
be a component of $\Omega-\Delta$, and let
$\Gamma_0={\rm stab}(\Omega_0,\Gamma)=
\{\gamma\in\Gamma;~\gamma(\Omega_0)=\Omega_0\}$, be the stabilizer of
$\Omega_0$ in $\Gamma$.
Since $\Gamma$ is a terminal regular torsion b-group of signature (2,0), we
have that $\Omega_0/\Gamma_0$ is an orbifold of signature
$(0,3;\infty,\infty,\infty)$. If $A_2^{1/2}(\Omega_0)\neq\Omega_0$, then
$\Gamma_0=\tilde{\Gamma}_0$, where
$\tilde{\Gamma}_0={\rm stab}(\Omega_0,\tilde{\Gamma})$,
and therefore the quotient orbifolds $S_0=\Omega_0/\Gamma_0$
and $\Omega_0/\tilde{\Gamma}_0$ are equal. If, to the contrary,
$A_2^{1/2}$ fixes $\Omega_0$, then it will induce a biholomorphic mapping in
$S_0$, say $f$. Since $A_2^{1/2}\not\in\Gamma$ but $A_2\in\Gamma$, we
have that $f$ is not trivial and has order $2$.
We also have that $f$ fixes the puncture determined by $A_2$, since
$A_2^{1/2}$ commutes with that element.
Therefore, $f$ has to interchange
the other two punctures and this implies that
$\Omega_0/<\Gamma_0,~A_2^{1/2}>=S_0/<f>$ has signature
$(0,3;\infty,\infty,2)$. Observe that we have used the fact that
$\tilde{\Gamma}_0$ is generated by $\Gamma_0$ and $A_2^{1/2}$.\end{pf}

The group $\tilde{\Gamma}$ has the following presentation
${\tilde{\Gamma}}=\{A_1,~A_3,~C_1,~C_3,~A_2^{1/2};~A_1,~A_2^{1/2},$\\
$A_3~$ are accidental parabolics,$~A_2^{-1/2}C_1^{-1},
C_1A_2^{1/2}A_1,
{}~A_1^{-1}A_2^{-1/2},~A_2^{1/2}A_3,~C_3A_2^{-1/2},\\
{}~A_2^{-1/2}C_3^{-1}A_3^{-1}$
are elliptic elements of order$\}$.
We will write only the expressions of the generators of
$\tilde{\Gamma}$ that we will use in this proof:
$$(C_1A_2^{1/2})=i\left[\begin{array}{cc}-1&2+\tau_1\\0&1\end{array}\right],~
A_1^{-1}=\left[\begin{array}{cc}-1&2\\0&-1\end{array}\right],~
A_2^{1/2}=\left[\begin{array}{cc}0&1\\-1&2\end{array}\right],$$
$$A_3^{-1}=\left[\begin{array}{cc}-1+2\tau_2(1-\tau_2)&2(1-\tau_2)^2\\
-2\tau_2^2&-1-2\tau_2(1-\tau_2)\end{array}\right],$$
$$C_3A_2^{-1/2}=i\left[\begin{array}{cc}-1+2\tau_2-\tau_2\tau_3+\tau_2^2\tau_3&
2-\tau_3-2\tau_2-\tau_2^2\tau_3+2\tau_2\tau_3\\
2\tau_2+\tau_2^2\tau_3&
1-2\tau_2+\tau_2\tau_3-\tau_2^2\tau_3\end{array}\right].$$

\sn {\bf 4.3.} Let $\cal F$ be a terminal regular b-group with signature
$(0,6;2,2,2,2,2,2)$
constructed by the techniques of section $3$, and corresponding
to the orbifold shown in figure $4$.

\noindent {\bf Remark}. It is easy to see that if
we apply yhe hyperelliptic
involution to the surface of the figure $3$, then we obtain a
partition in the quotient surface as given in figure $4$. That is why
we have chosen this partition among all the possible ones on an
orbifold with signature $(0,6;2,\ldots,2)$.

A set of generators for $\cal F$ (equivalent to the generators of
$\tilde{\Gamma}$ above
computed) consists of the following transformations:
$$D_1=i\left[\begin{array}{cc}-1&0\\0&1\end{array}\right],~
B_1=\left[\begin{array}{cc}-1&-2\\0&-1\end{array}\right],~
B_2=\left[\begin{array}{cc}-1-\alpha&\alpha^2\\-1&-1+\alpha
\end{array}\right],$$
$$B_3=\left[\begin{array}{cc}-1+2\beta+2\alpha\beta^2&
-2(1+\alpha\beta)^2\\
2\beta^2&-1-2\beta-2\alpha\beta^2
\end{array}\right],$$
$$D_4=i\left[\begin{array}{cc}-1-2\alpha\beta+2\alpha\gamma+2\gamma\beta^{-1}&
-2(1+\alpha\beta)(-\alpha\beta^2+\gamma+\alpha\beta\gamma)\beta^{-2}\\
2\gamma-2\beta&
1+2\alpha\beta-2\alpha\gamma-2\gamma\beta^{-1}
\end{array}\right].$$
$\alpha,~\beta$ and $\gamma$ are three complex numbers, chosen so that the
above matrices have nice expressions.

To complete the proof of the theorem, we have to find a \M transformation $E$
such that $E{\cal F}E^{-1}=\tilde{\Gamma}$,
and then we have to express the coordinates
of $\cal F$ in terms of the numbers $\alpha~,\beta$ and $\gamma$.
Topological considerations give that $EB_1E^{-1}=A_1$ and
$EB_2E^{-1}=A_2^{-1/2}$.
This implies that $E(z)=-z+1+\alpha.$
We also have $ED_1E^{-1}=A_2^{-1/2}C_1^{-1}$ which gives $\alpha=\tau_1/2$.
It is a matter of computation to see that the conjugation
$EB_3E^{-1}=A_3^{-1}$ gives the following four equations,
whose unique answer is
$\beta=\tau_2$:
$$\left\{\begin{array}{lcl}
-1+2\beta+2\alpha\beta^2+2\beta^2b&=&1+2\tau_2(1-\tau_2),\\
-2(1+\alpha\beta)^2+b(-1-2\beta-2\alpha\beta^2)&=&
-b(-1-2\tau_2+2\tau_2^2)-2(1-\tau_2)^2,\\
-2\beta^2&=&-2\tau_2^2,\\
1+2\beta+2\alpha\beta^2&=&-2\tau_2^2b-1+2\tau_2(1-\tau_2).
\end{array}\right .$$
The relation $ED_4E^{-1}=C_3A_2^{-1/2}$ gives the set of equations
$$\left\{\begin{array}{lcl}
-1-2\alpha\beta+2\alpha\gamma+2\gamma\beta^{-1}+b(2\gamma-2\beta)&=&
-1+2\tau_2-\tau_2\tau_3+\tau_2^2\tau_3,\\
2\beta-2\gamma&=&2\tau_2+\tau_2^2\tau_3.\end{array}\right . ,$$
whose answer is
$\gamma={\load{\large}{\displaystyle}-\frac{\tau_2^2\tau_3}{2}}$.

In \cite{kra:horoc} is proven that $(\tau_1,\tau_2,\tau_3)$
is a set of coordinates
for $T(2,0)$. The coordinates for $T(0,6;2,2,2,2,2,2)$ are given by
$z_1=\alpha$, $z_2=1+\beta$ and $z_3=1-\frac{\gamma}{\beta^2}$. Substituting
the values of $\alpha$, $\beta$ and $\gamma$ obtained above we get
that the expression of the isomorphism is the one given in the
statement of theorem $11$.~\hfill$\Box$

\sn {\bf 4.4.} As a corollary of the proof of the theorem $11$
we obtain the other isomorphisms of \S $4.1$.
\begin{cor}The mappings $$\tau_1\mapsto \frac{\tau_1}{2}~~
{\rm and}~~
(\tau_1,\tau_2)\mapsto (\frac{\tau_1}{2},1+\tau_2),$$
give the isomorphisms $T(1,1;\infty)\cong
T(0,4;\infty,2,2,2)$ and $T(1,2;\infty,\infty)\cong
T(0,5;\infty,2,\newline\ldots,2)$, respectively, for some choice of horocyclic
coordinates.\end{cor}
\begin{pf}
The argument goes as follows: to construct the surface of genus 2 with the
partition given in the figure $3$, we
start with a thrice punctured sphere, $S_1$;
then we glue two of the punctures, obtaining a surface $S_2$
with signature $(1,1;\infty)$. This construction uses only the
coordinate $\tau_1$, and therefore it gives the first isomorphism of
the theorem. The next step is to
glue to the puncture of $S_2$ a three times punctured sphere
to get a surface $S_3$, with signature $(1,2;\infty,\infty)$.
For this construction we need the coordinates $(\tau_1,\tau_2)$. Thus
we obtain the second isomorphism, completing the proof of the
theorem.\end{pf}

This technique can be used to compute different isomorphisms between \T spaces.
For example, there is another partition of a surface of genus 2;
the hyperelliptic
involution can be found in \cite{kra:horoc}. With computations similar to the
ones described in this section, one can find the Patterson isomorphisms
in that case. By a result of Kravetz, the set of fixed points of
a biholomorphic map of finite
order on a \T space is isomorphic to another \T space of lower dimension (see
\cite{krav:teic} or \cite[pp. 259-260]{nag:teic}).
Some of those isomorphisms can also be studied with our techniques.

\ifx\undefined\bysame
\newcommand{\bysame}{\leavevmode\hbox to3em{\hrulefill}\,}
\fi

\sn School of Maths, Tata Institute of Fundamental Research,
Bombay, India\\
pablo@motive.math.tifr.res.in
\end{document}